\documentclass[12pt,draft]{elsarticle}

\usepackage{epic,eepic}
\usepackage{amssymb,amsfonts,tikz}
\usepackage{amsthm}
\usepackage{mathrsfs}
\usepackage{fancybox}
\usepackage{fancyhdr}
\usepackage{mathtools}

\journal{Linear Algebra and its Applications}

\newtheorem{thm1}{Theorem}
\newtheorem{corollary}[thm1]{Corollary}
\newtheorem{proposition}[thm1]{Proposition}
\newtheorem{remark}[thm1]{Remark}

\newtheorem{example}{Example}
\newtheorem{lemma}[thm1]{Lemma}

\begin{document}
	
	\begin{frontmatter}
		
		
		
		\title{Spectra of Tridiagonal Matrices over a Field}
		
		
		\author{R. S. Costas-Santos\corref{cor1}}
		\ead{rscosa@gmail.com}
		\ead[url]{http://www.rscosan.com}
		\cortext[cor1]{Corresponding author}
		\address{Dpto. de F\'isica y Matem\'aticas, 
			Facultad de Ciencias, Universidad de Alcal\'a, 
			28871 Alcal\'a de Henares, Spain}
		
		\author{C. R. Johnson}
		\ead{crjohnso@math.wm.edu}
		\address{Department of Mathematics,
			College of William and Mary, Williamsburg, 
			VA 23187}
		

\bibliographystyle{plain}

\pagestyle{myheadings}
\markboth{R. S. Costas-Santos et al.}
{Spectra of Tridiagonal Matrices over a Field}

\begin{abstract}
We consider spectra of $n$-by-$n$ irreducible tridiagonal 
matrices over a field and of their $n-1$-by-$n-1$ 
trailing principal submatrices.  The real symmetric and 
complex Hermitian cases have been fully understood: it 
is necessary and sufficient that the necessarily real 
eigenvalues are distinct and those of the principal 
submatrix strictly interlace. So this case	is very restrictive.

By contrast, for a general field, the requirements on the 
two spectra  are much less restrictive. In particular, in 
the real or complex case, the $n$-by-$n$ characteristic 
polynomial is arbitrary (so that the algebraic multiplicities 
may be anything in place of all 1's in the classical cases) 
and that of the principal submatrix is the complement of a 
lower dimensional algebraic set (and so relatively free). 
Explicit conditions are given. 
\end{abstract}

\begin{keyword}
Eigenvalues\sep Irreducible\sep Orthogonal polynomials 
\sep Characteristic polynomial\sep  Recurrence 
relation\sep Tridiagonal matrix.
\MSC[2010] Primary 15A18 \sep 15B05 \sep 05C05\sep 42C05.
\end{keyword}

\end{frontmatter}


\begin{abstract}

\end{abstract}
\section{Introduction} \label{intro-sec}
An $n$-by-$n$ matrix $A=(a_{ij})$ is called {\tt tridiagonal} 
if $|i-j|>1$ implies $a_{ij}=0$. Such a matrix may have 
nonzero entries only on the sub-, super-, and main diagonals.

\begin{center}
\begin{tikzpicture}[domain=-2:2,scale=0.7,samples=200]
\draw (0,0) node {$A=$}
(0.7,0)++(0.2,2)-- ++(-0.2,0)-- ++(0,-4)-- ++(0.2,0)
(5,0)++(-0.2,2)-- ++(0.2,0)-- ++(0,-4)-- ++(-0.2,0)
(1,1.7) -- (4.7,-1.7) (2,1.7) -- (4.7,-0.7)  (1,0.7) -- (3.7,-1.7)
(4,1.1) node {\Large $0$}   (1.7,-1.1) node {\Large $0$}
 ;
\end{tikzpicture}
\end{center}

We are interested in the eigenvalues of such a matrix 
and of its trailing $(n-1)$-by-$(n-1)$ principal submatrix: 
$A(1)=A[\{2,3, \dots, n\}]$,  which we view in terms of their 
characteristic polynomials, over a general field $\mathbb F$. 
If $\{\lambda_1, \cdots, \lambda_n\}$ occur as the 
eigenvalues of $A$ and $\{\mu_1, \cdots, \mu_{n-1}\}$ 
as the eigenvalues of $A(1)$, they will also occur for a tridiagonal
matrix with all super-diagonal entries nonzero. 
When the super-diagonal entries are all nonzero, they may be 
normalized to be all be 1's, via diagonal similarity.

So, wlog, we consider {\tt normalized tridiagonal} matrices.
Our $A$ looks like 
\begin{equation} \label{MatA} 
A=\left[\begin{array}{cccccc} a_{1} & 1 &0 & \cdots &  0 \\ 
b_{1} & a_{2} & 1 & \ddots & \vdots \\ 
0 & b_{2} & \ddots & \ddots & 0\\ 
\vdots & \ddots & \ddots & \ddots & 1 \\ 
0 & \cdots & 0 & b_{n-1} & a_{n}
\end{array}\right]
\end{equation} 

In this  case, some sub-diagonal entries may be 0, in which case $A$ is 
{\tt reducible}, or all may be nonzero, in which case $A$ is {\tt irreducible}.
In the reducible case, $A$ and $A(1)$ must have eigenvalues in common. 

The number of common eigenvalues is $k$ (counting multiplicities)
if and only if $b_{n-k}=0$ and no prior $b_i's$ are 0.  We focus upon 
the irreducible case, from which the reducible case may be deduced,
in which case $A$ and $A(1)$ have no common eigenvalues. 

Let $p_n(t)=\det (t I-A)$ and $p_{n-1}(t)=\det (tI - A(1))$, the 
characteristic polynomial of $A$ and $A(1)$, respectively. 
Generally let $p_k(t)$ be the characteristic polynomial of the trailing 
$k$-by-$k$ principal submatrix.

Via determinantal expansion, we have the following known relationships:  
\begin{equation} \label{1:1} 
p_n(t)=(t-a_1)p_{n-1}(t)-b_1 p_{n-2}(t),
\end{equation} 
and generally
\begin{equation} \label{1:2} 
p_{k+1}(t)=(t-a_{n-k})p_{k}(t)-b_{n-k} p_{k-1}(t),\quad 
k=0, \dots, n-1,
\end{equation} 
in which $p_{-1}(t)=0$ and $p_{0}(t)=1$.

From these, it is clear, in the irreducible case, that 
\begin{enumerate}
\item[i)] $p_{n-2}$, $p_{n-3}$, $\dots$, $p_1$ are uniquely 
determined by $p_n$ and $p_{n-1}$, and thus the eigenvalues 
of all the trailing principal submatrices are determined by 
those of the first two; 
\item[ii)] $p_n$ and $p_{n-1}$ are relatively prime, so 
that $A$ and $A(1)$ have no eigenvalues in common. The 
same is true for $p_{k+1}$ and $p_k$. However, $p_{k+1}$ 
and $p_{k-r}$ could have  common roots for $r\ge 1$.
\begin{example}
The matrix 
\[
\left[\begin{array}{ccc}
0 & 1 & 0 \\ 1 & 0 & 1 \\ 0 & 1 & 0
\end{array}\right]
\]
has spectra $\sigma(A)=\left\{-\sqrt 2,0,\sqrt 2\right\}$ and 
$\sigma(A(1,2))=\{0\}$.
\end{example}
\item[iii)] Unlike the Hermitian case, $p_n$ may have 
multiple roots but, like the Hermitian case, each of them 
has geometric multiplicity 1 as an eigenvalue of $A$ (as 
rank$(tI-A)\ge n-1$ for each $t$).

So the eigenvalues of $A$ may be algebraically multiple, 
but not geometrically so.
\end{enumerate} 

Suppose now that $p$ and $q$ are given monic polynomials, 
i.e., the leading coefficient is 1, over $\mathbb F$, of 
degree $n$ and $n-1$, respectively. If there is an 
$n$-by-$n$ tridiagonal matrix over $\mathbb F$, 
such that $p=p_n$ and $q=p_{n-1}$, we call $p$ and 
$q$ a {\tt tridiagonal pair} (TrP), and if the 
tridiagonal matrix may be taken to be irreducible, we 
call $p$ and $q$ an  {\tt irreducible tridiagonal 
pair} (ITrP).

We seek to understand which pairs are ITrP over $\mathbb F$ 
and which polynomials $p$ occur in an ITrP. Not all pairs 
are ITrP, but every monic $p$ of degree $n$ does occur 
as $p_n$ when $\mathbb F=\mathbb R$ or $\mathbb C$ (not 
in general). 
So, over $\mathbb R$ or $\mathbb C$ an irreducible 
tridiagonal matrix may have any characteristic polynomial, 
and thus, any algebraic multiplicities for its eigenvalues,  
in stark contrast to the real symmetric or Hermitian cases. 
\section{Theorems and Examples}
Given monic polynomial $p$ and $q$ over $\mathbb F$ of 
degrees $n$ and $n-1$, when we apply the division 
algorithm to them, one of two things may happen:
either a) the degree of the remainder drops by 
exactly 1 each time, so that the algorithm 
consumes $n-1$ steps, or b) at some stage of 
the division algorithm, there is a drop in 
degree by more than 1. 

In case a), which is generic over $\mathbb R$ or 
$\mathbb C$, we call $p$ and $q$ a {\tt proper pair}. 
If $p$ is given, we say that $q$ {\tt is proper with} 
$p$, and we may refer to the set of such $q$  as the 
{\tt proper set for} $p$. In case b), an algebraic 
condition must be satisfied by the coefficients of $p$ 
and $q$. There are no more than $n-2$ of these.
Thus, the non-proper pairs form an algebraic set, 
and the proper set of $p$ is also algebraic. It 
follows that the proper pairs are the complement of 
an algebraic set, and likewise for the proper set 
of $p$.

We may now observe a basic characterization.
\begin{thm1} \label{thm:1} 
Let $p$ and $q$ be monic polynomials, over a field 
$\mathbb F$, of degree $n$ and $n-1$ respectively. 
Then $p$ and $q$ form an ITrP if and only if they 
are a proper pair.
\end{thm1} 

\begin{proof} 
Suppose that $p$ and $q$ form a proper pair. Then, 
upon division of $p$ by $q$, according to (\ref{1:1}), 
we may conclude  what $a_1$, $b_1$ and $p_{n-2}$ would 
have to be in order to have $p_n=p$ and $p_{n-1}=q$. 

Since $p$, $q$ is proper, $b_1$ is nonzero and well-
defined, and $\deg p_{n-2}=n-2$. However, also since 
$p$, $q$ is proper, we may continue by applying 
(\ref{1:2}) to $p_{n-1}=q$ and $p_{n-2}$ to get 
$a_2$, $b_2$ and $p_{n-3}$ in the same way; $b_2\ne 0$ 
and $\deg p_{n-3}=n-3$. Again, as $p$, $q$ is proper, 
we may continue to get $a_3$, $b_3$ and $p_{n-4}$, 
and so on. 

This allows us to construct the unique (normalized) 
irreducible tridiagonal matrix $A$ for which $p_{n}=p$ 
and $p_{n-1}=q$, which shows that $p$, $q$ is an ITrP.

If $p$ and $q$ form an ITrP, the proof that $p$ and 
$q$ is a proper pair is similar. We have $p_n=p$ and 
$p_{n-1}=q$, so that (\ref{1:1}) and (\ref{1:2}) 
imply that $p$, $q$ is a proper pair, as $b_1$, $b_2$, 
$\dots$, $b_{n-1}\ne 0$.
\end{proof} 

We note that when $p$, $q$ is a proper pair (and thus 
an ITrP), the irreducible (normalized) tridiagonal 
matrix that realizes them is uniquely determined. 
So $p_{n-2}$, $\dots$, $p_1$ (and their roots, the 
eigenvalues  of the trailing principal submatrices) are 
fully determined.  
It is an interesting question how these roots are a 
function of the roots of $p$ and $q$.

We note that not every relatively prime pair $p$ and 
$q$ is an ITrP, even over $\mathbb R$.

\begin{example} 
Let $n=3$ and $2$, $-3$, $-5$ be the 
roots of monic $p_3$ and $1$, $-1$ the roots of the
monic $p_2$. Then $p_2$ and $p_3$ are relatively prime, 
but there is no tridiagonal matrix with eigenvalues 
$2$, $-3$, $-5$ and with $1$, $-1$ as the eigenvalues 
of the upper let $2-$by$-2$ principal  submatrix. 
We have 
\[
p_3(t)=t^3+6t^2-t-30, \quad and \quad 
p_2(t)=t^2-1.
\]
Suppose there is a tridiagonal matrix $A=(a_{ij})$. 
Then
\[
p_3(t)=(t-a_{1})p_2(t)-b_{1}p_1(t).
\]
This implies that $a_{1}=-6$ and that $p_1(t)=24/b_{1}$, 
which is  a  polynomial of degree 0. Therefore 
$p_3$, $p_2$ is not a proper pair and not an ITrP.
\end{example}

\begin{example}
Let $p_4$ have roots $-1$, $-2$, $3$, $4$ and 
$p_3$ have roots $-3$, $1$, $2$. Then $p_3$ and 
$p_4$ are relatively prime. Let us assume that 
there is a $4$-by-$4$ tridiagonal matrix 
$A$ with  $p_A(t)=p_4(t)$ and $p_{A(1)}(t)=p_3(t)$.
If we apply the division algorithm to  $p_3$ 
and $p_4$, we get 
\[
p_4(t)=(t-4)p_3(t)-12t+48=(t-4)p_3(t)-12p_2(t).
\]
Therefore the degree of $p_2$ drops by 2; hence 
such a $4-$by$-4$ tridiagonal 
matrix does not exist.
\end{example}

However, because the proper set of a monic polynomial 
over $\mathbb R$, or $\mathbb C$ is the complement of 
a sufficiently low dimensional algebraic set, the 
proper set is necessarily nonempty.

\begin{thm1} \label{thm:2} Suppose that $p$ is a 
monic degree $n$ polynomial over $\mathbb R$ or 
$\mathbb C$. Then, there is a monic polynomial $q$ 
over the same field as $p$ such that $p$ and $q$ 
form an ITrP.
\end{thm1} 

\begin{proof}
For a given $p$ the existence of such $q$ is 
sraightforward since the proper set of $p$ is the 
complement an algebraic set which is strictly 
contained in $\mathbb R^{n-1}$ or $\mathbb C^{n-1}$. 
\end{proof}	

Nevertheless, it may happen for other fields that 
the result for theorem \ref{thm:2} is not true.

\begin{example}
Over the field $GF_2$, not every monic polynomial 
is attained as the characteristic polynomial of an 
irreducible $3$-by-$3$ tridiagonal matrix.

Let $p(t)=t^3+1$.  If $A\in M_3(GF_2)$ is 
irreducible and tridiagonal, then $A$ is of the form 
\[
A=\left[ \begin{array}{ccc}d_1 & 1 & 0 \\1 & d_2 
& 1 \\ 0 & 1 & d_3 \end{array} \right],
\]
with each $d_i$ =0 or 1, $i=1, 2, 3$. Then 
\[
p_A(t)=t^3-(d_1+d_2+d_3)t^2+(d_1d_2+d_1d_3+d_2d_3)t-
(d_1d_2d_3-d_1-d_3).
\]
For $d_1+d_2+d_3=0$, either 0 or 2 of the $d_i's$ 
is 1. Then for $d_1d_2+d_1d_3+d_2d_3=0$, it must 
be that all $d_i=0$ (if two are equal to 1, this 
expression is 1). But if all are 0, then 
$d_1d_2d_3-d_1-d_3=0$, not 1.
\end{example}

By a simple counting argument, over any finite field 
$\mathbb F$ some polynomials do not occur as the 
characteristic polynomial of a normalized irreducible 
tridiagonal matrix. If $\mathbb F$ has $k$ 
elements, then there are $(k-1)^{n-1}$ such matrices,
but $k^n$ distinct monic polynomials. It is an interesting 
question which polynomials are realized.

\begin{corollary} \label{cor:1}
Over $\mathbb R$ or $\mathbb C$, an irreducible 
tridiagonal matrix may have any characteristic 
polynomial (and thus, any eigenvalues, 
counting multiplicities).
\end{corollary}

For real symmetric and complex Hermitian 
irreducible tridiagonal matrices, it is known 
\cite{Hoc} that the only multiplicity list that 
occurs for the eigenvalues is all 1's.
And, in general, the maximum geometric 
multiplicity is 1. However, for algebraic 
multiplicity, the situation is quite different.

For a further reading about the 
Inverse eigenvalue problems
for band matrices see e.g. \cite{bogo, bogo1, hijo}.

\begin{corollary} \label{cor:2}
Any partition of $n$ may be the list of algebraic 
multiplicities of  an irreducible tridiagonal 
matrix over $\mathbb R$ or $\mathbb C$.
\end{corollary}

The (undirected) graph of an irreducible tridiagonal 
matrix is simply a path. We conjecture that the 
same is true for other trees, i.e. any algebraic 
multiplicities may occur, and this is true 
for the star on $n$ vertices \cite{Hil, JoLe}.

Though over $\mathbb R$ or $\mathbb C$ any
polynomial occurs as the characteristic 
polynomial of an irreducible tridiagonal matrix, it is 
not easy to explicitly give a tridiagonal matrix realization.
In the next section, we show how a realization may 
be given, using some ideas form orthogonal 
polynomials.
\section{Tridiagonal Matrices and Orthogonal Polynomials}
The theory of linear functionals is a natural tool 
to understand tridiagonal realizability. 
We first give some basic facts we need.

Given a linear functional ${\mathscr L}:\mathbb F[t] \to 
\mathbb F$ we denote by $m_k={\mathscr L}(t^k)$, for 
all $k=0, 1, \dots$, the {\tt moments} of $\mathscr L$, and 
by $H_k$ the  $(k+1)$-by-$(k+1)$ Hankel matrix
\begin{equation} \label{3:3}
H_k=\left[\begin{array}{lllll} 
m_0 & m_1 & m_2 & \cdots & m_k \\
m_1 & m_2 & m_3 & \cdots & m_{k+1} \\ 
\vdots &\vdots &\vdots &\ddots &\vdots \\
m_k & m_{k+1} & m_{k+2} & \cdots & m_{2k} \\ 
\end{array}\right].
\end{equation} 
The linear functional ${\mathscr L}$ is said to be 
{\tt quasi-definite} if $\det(H_k)\ne 0$ for 
all $k=0, 1, \dots$. 

\begin{remark}
In this work we fix an integer, $n>0$, and, since we 
are interested in  $n$-by-$n$ matrices, it is enough 
to suppose that $\det H_k\ne 0$ for $k=0, 1, \dots n-1$, 
Therefore it is not an issue if there exists  some $N>n$ 
such that $H_N$ is singular.
So, we will say that $\mathscr L$ is quasi-definite 
if the matrices $H_0$, $H_1$, $\dots$, $H_{n-1}$ are 
all invertible.
\end{remark}

The following result is well-known for orthogonal 
polynomial sequences:
\begin{proposition} \cite[p. 17]{chi}
For any quasi-definite linear functional ${\mathscr L}$, 
there exists a polynomial sequence $\{p_k\}$, unique up 
to a multiplicative constant,  defined by 
\begin{equation} \label{3:4}
p_k(t)=\left[\begin{array}{lllll} 
m_0 & m_1 & m_2 & \cdots & m_k \\ m_1 & m_2 & m_3 
& \cdots & m_{k+1} \\ \vdots &\vdots &\vdots &\ddots 
&\vdots \\ m_{k-1} & m_{k} & m_{k+1} & \cdots 
& m_{2k-1} \\ 1 & t & t^2 & \cdots & t^{k} 
\end{array}\right],\quad k=0, 1, \dots,
\end{equation} 
that fulfills the property of orthogonality
\[
{\mathscr L}(p_\ell p_k)=0,\quad n\ne m,\quad 
\ell,k  =0, 1, \dots,
\]
\[
{\mathscr L}(p_k^2)\ne 0,\quad k =0, 1, \dots
\]
\end{proposition}
Note that, due to the normalization taken for 
our tridiagonal matrices, we need to consider
the following normalization for the polynomials:
\[
P_0(t)=1,\quad P_k(t)=\frac 1{\det(H_{k-1})} p_k(t),
\quad k=1, 2, 3, \dots.
\]
and $P_k(0)=\det(\widetilde{H}_{k-1})/\det(H_{k-1})$, 
in which
\[
\widetilde{H}_{k-1}=\left[\begin{array}{llll} 
m_1 & m_2 & \cdots & m_{k} \\
m_2 & m_3 & \cdots & m_{k+1} \\ 
\vdots &\vdots  &\ddots &\vdots \\
m_k & m_{k+1}  & \cdots & m_{2k-1} \\ 
\end{array}\right].
\]
\begin{remark}
Since we are considering tridiagonal matrices normalized
so that the superdiagonal is all 1's, it is more 
convenient to use monic polynomials. Observe that 
if $(p_k)$ satisfies the recurrence relation 
(\ref{1:2}) then $(P_k)$ satisfies the following 
recurrence relation for $k=1, \dots, n-1$:
\begin{equation}\label{3:5}
tP_k(t)= \frac{\det (H_k)}{\det (H_{k-1})}P_{k+1}(t)+
a_{n-k} P_k(t)+b_{n-k} \frac{\det (H_{k-2})}
{\det (H_{k-1})} P_{k-1}(t),
\end{equation}
and since ${\mathscr L}$ is quasi-definite, 
$\det(H_k)\ne 0$ for all $k=0, 1, \dots, n-1$.

Observe that $p_k(x)$ has degree $k$ if and only if 
$H_k$ is regular for $k=0, 1, \dots, n-1$. Hence, the 
fact that $\mathscr L$ is quasi-definite 
means that $p_{n}$ and $p_{n-1}$ form a
proper pair.
\end{remark}

For further reading on the existence of orthogonal 
polynomial sequences and this matrix representation 
we suggest \cite[Charper 3]{chi}.

Taking all this into account now we can state an 
explicit result about when $p$ and $q$ form an ITrP.

\begin{thm1} \label{thm:3}
For any polynomials $p$ and $q$ of degree $n$ and $n-1$,
respectively, with coefficients over a field $\mathbb F$, let us 
denote $p$ by $p_n$ and $q$ by $p_{n-1}$, and let us 
consider the following two linear functionals:
\begin{itemize} 
\item If all the roots of $p$ are different, then
${\mathscr L}_1: \mathbb F[t]\to \mathbb F$ 
\[
{\mathscr L}_1(f(t))=\sum_{k=1}^n \frac{f(\lambda_k)}
{p'_n(\lambda_k)p_{n-1}(\lambda_k)},
\]
where
\[
p_n(t)=\prod_{k=1}^n (t-\lambda_i), \quad and \quad 
p_{n-1}(t)=\prod_{k=1}^{n-1}(t-\mu_i).
\]	 
\item If all the roots of $p$ are the same, namely 
$a$ with multiplicity  $n$, then 
${\mathscr L}_2: \mathbb F[t]\to \mathbb F$ 
\begin{equation} \label{3:7}
{\mathscr L}_2(f(t))=C \frac{d^{n-1}}{z^{n-1}}
\left(\frac{f(z)} {p_{n-1}(z)}\right)(a),
\end{equation} 
where $C$ is a constant such that ${\mathscr L}(1)=m_0$.
\end{itemize}

Then, the following statements are equivalent:
\begin{enumerate}
\item $p$ and $q$ form a proper pair.
\item All the Hankel matrices $H_k$ associated with 
the linear functional $\mathscr L$ are invertible 
for $k=0, 1, \dots, n-1$.
\end{enumerate} 
\end{thm1}

\begin{proof}
WLOG we need to prove this result in the following 
two situations:
i) when all the zeros of $p$ are different, and 
ii) when $p$ has one zero with multiplicity $n$.

Let $\lambda_{1},  \dots, \lambda_{n}$ be the zeros 
of $p$, all of them different, and $\mu_{1},\dots, 
\mu_{n-1}$ be the zeros of $q$ over the field 
$\mathbb F$, such that 
$\{\lambda_1, \cdots, \lambda_n\} \cap 
\{\mu_1,\cdots, \mu_{n-1}\}=\emptyset$. 

By definition, the functional $\mathscr L_1$ is linear. 
Moreover, if it is quasi-definite then the Hankel 
matrices associated to it are invertible, i.e. 
$\det(H_k)\ne 0$ for $k=0, 1, \dots, n-1$. 

Let us define the $n$-by-$n$ tridiagonal matrix $A$ 
defined in \eqref{MatA} where 
\[
b_{i}=\frac{{\mathscr L}_1(p_{n-i}^2)}
{{\mathscr L}_1(p_{n-i-1}^2)}, \qquad  i=1, 2, \dots, n-1,
\]
and 
\[
a_{i}= 
\frac{-p_{n-i+1}(0)-b_{i} p_{n-i-1}(0)}{p_{n-i}(0)}, 
\qquad  i=1, 2, \dots, n,
\]
being $p_{-1}(t)=0$ and $p_0(t)=1$.

\begin{remark}
Note that we can assume $p_k(0)\ne 0$ for all $k$, because 
if not we apply a linear change of variables $y(x)=x+b$, 
$b\ne 0$, so that the recurrence relation coefficients $b_i's$ 
remain the same and, since $p_k(b)\ne 0$ for all 
$k=1,\dots,n-1$, then 
\[
a_{i}= 
\frac{-p_{n-i+1}(b)-b_{i} p_{n-i-1}(b)}{p_{n-i}(b)}
\]
are finite.
\end{remark} 
If we prove that $p_A(t)=p_n(t)$ and 
$p_{A(1)}(t)=p_{n-1}(t)$
and the matrix is irreducible, i.e., $b_{i}\ne 0$ 
for all $i=1, \dots, n-1$, then the necessary condition 
holds and, therefore, $p_n(t)$, $p_{n-1}(t)$ form a ITrP. 

By construction we know there exists $(q_k(t))_{k=0}^n$ 
a sequence
of monic polynomials orthogonal with respect to $\mathscr 
L_1$, i.e., they fullfills the following property of 
orthogonality:
\[
{\mathscr L}_1 (t^\ell q_k(t))=0,\qquad \ell=0, 
1, \dots, k-1, \qquad k=1, 2, \dots, n,
\]
as well as the three-term recurrence relation
\begin{equation} \label{3:6} 
t q_k(t)=q_{k+1}(t)+\alpha_{n-k} q_k(t)+\beta_{n-k} 
q_{k-1}(t), \quad k=1, 2, \dots, n-1.
\end{equation} 
By using the previous recurrence relation and the 
orthogonality conditions for $q_k$, it is 
straightforward to prove $q_n(t)=p_n(t)$, as well 
as $\beta_1=b_1$. Moreover, if we set $t=0$ and 
$k\mapsto n-1$ in (\ref{3:6}) we obtain that 
$\alpha_{1}=a_{1}$. 
In order to prove that  $q_{n-1}(t)=p_{n-1}(t)$ we 
need the followig result.
\begin{lemma}
For any polynomial $p(t)$ of degree $m>1$, with 
different zeros $x_1, x_2, \cdots, 
x_m$, the following identity holds true:
\[
\sum_{j=1}^{m} \frac 1{p'(x_j)}=0.
\]
\end{lemma}
So 
\[
{\mathscr L}_1(p_{n-1}(t))=\sum_{j=1}^n\frac{1}
{p_n'(\lambda_j)}=0,
\]
and if we consider, for $\ell=1,\dots, n-2$, the polynomials 
$\pi_\ell(t)=(t-\lambda_1) \cdots (t-\lambda_\ell)$,  then 
\[
{\mathscr L}_1(p_{n-1}(t)\pi_\ell(t))=
\sum_{j=\ell+1}^n \frac{\pi_\ell(\lambda_j)}{p_n'(\lambda_j)}=
\sum_{j=\ell+1}^n \frac{1}{(p_n/\pi_\ell)'(\lambda_j)}=0,
\]
therefore, by unicity, $q_{n-1}(t)=p_{n-1}(t)$.

In fact, since
\[
{\mathscr L}_1(q^2_k(t))={\mathscr L}_1(t^kq_k(t))=\beta_k 
{\mathscr L}_1(t^{k-1}q_{k-1}(t))=
{\mathscr L}_1(1) \beta_1\beta_2\cdots \beta_k\ne 0,
\]
we get that, by construction, for $k=1, 2, \dots, n$, 
$\beta_{n-k}=b_{n-k}$ and by the orthogonlity conditions 
$\alpha_{n-k}=a_{n-k}$. Therefore $p_A(t)=p(t)$ and 
$p_{A(1)}(t)=q(t)$.

Remember that, by construction, we have
\[
0\ne \big(\det(H_{k-1})\big)^2{\mathscr L}_1(P^2_k)=
{\mathscr L}_1(p^2_k)= 
b_k b_{k-1}\cdots b_1  {\mathscr L}_1(1).
\]

And it is sufficient to have a tridiagonal pair, 
because in such a case there exists a matrix $A$ 
so that $p_A(t)=p_n(t)$ and $p_{A(n)}(t)=p_{n-1}(t)$. 
So we consider the same inner product and, by 
construction, the polynomial $p_k(t)$ has degree 
$k$ for $k=0, 1, \dots, n$, and they are monic. 

Then if we establish the orthogonality conditions 
again, we check in a straightforward way that the 
leading coefficient of the matrix expression (\ref{3:4})
is indeed $\deg(H_k)$ that must be non-zero so, 
the linear functional is quasi-definite and that completes 
the proof for this case.

If $p$ has one zero, namely $a$, with multiplicity 
$n$, then we consider the the linear functional 
${\mathscr L}_2$. 

Since the key to the proof is not about the expression 
for $\mathscr L_2$ but about the fact that the operator 
is linear we leave this part of the proof to the reader.
\end{proof}

\begin{example} \label{ex:5}
If $p(x)=x(x-1)^3(x+2)(x-5)$, we need to consider the
linear functional that is a linear combination of the ones
presented in theorem \ref{thm:3}, i.e.
\[
{\mathscr L}(f)=\frac{f(0)}{10}-\frac{f(-2)}{378}
+\frac{f(5)}{80640}
+\frac 1{\pi i}\int_{|z-1|=1} \frac{f(z)}{(z-1)^3 (z+1)^2}\, dz,
\]
where we have considered for the construction of the 
coefficients of the first part the polynomial $(x+1)^2$, 
but any polynomial of degree 2 or greater, proper with 
$x(x+2)(x-5)$, can be chosen.

With this construction we get the following sequence of 
polynomials:
\begin{eqnarray*}
p_0(x)&=&1 \\[2mm]
p_1(x)&=&x+\frac{811}{2193}
\\[2mm]
p_2(x)&=& x^2+\frac{4310 x}{2199}+\frac{7097}{6597}
\\[2mm]
p_3(x)&=&x^3-\frac{121347 x^2}{39845}
+\frac{26393 x}{7969}-\frac{37047}{39845}
\\[2mm]
p_4(x)&=&x^4-\frac{73660 x^3}{33301}+\frac{11933 x^2}
{33301}+\frac{33874 x}{33301}-\frac{40680}{33301}\\[2mm]
 p_5(x)&=&x^5-\frac{50023 x^4}{8243}+\frac{14679 x^3}{8243}
 +\frac{175435 x^2}{8243}-\frac{236462 x}{8243}+\frac{93312}{8243}
 \\[2mm]
p_6(x)&=&x(x-1)^3(x+2)(x-5)=p(x).
\end{eqnarray*}
Observe that, by construction, $q(x)=p_5(x)$ is proper with 
$p(x)$. 
Moreover, we obtain the 6-by-6 tridiagonal matrix 
\eqref{MatA} where 
\begin{eqnarray*}
\vec a=&\hspace{-3mm}\left(-\dfrac{565}{8243},
\dfrac{1058636543}{274500143},-\dfrac{1105993747}{1326878345},
\dfrac{438574003}{87619155},-\dfrac{852049}{535823},
-\dfrac{811}{2193}\right) \\[3mm]
\vec b=&\hspace{-3mm} \left(-\dfrac{43158096}{67947049},\dfrac{7882616040}{1108956601},\dfrac{659060091}{1587624025},-\dfrac{58253390}{4835
   601},-\dfrac{2345600}{4809249} \right).
\end{eqnarray*}
\end{example}
\begin{remark} Note that in the proper case there is 
an iterative algorithm to construct the realizing 
tridiagonal matrix computationally.
\end{remark} 

\begin{example} \label{ex:6} 
Here, we want  to give an example in which $p$ has 
degree 3 and multiple roots, and  $q$ is of degree 2. 
We obtain conditions for them to be proper pair can.

Consider the polynomials $p(x)=(x+1)(x-1)^2$, 
and $q(x)=(x-a)(x-b)$, $a, b \ne \pm 1$. With these 
polynomials we define the linear functional 
\[
{\mathscr L}_{a,b}(f)=\frac{f(-1)}{4(1+a)(1+b)}+\frac{f'(1)}{1-b}
-+\frac{f(1)}{(1-b)^2}.
\]
After a straightforward calculation we get the determinant 
of Hankel matrices for this linear functional
\begin{eqnarray*}
\det(H_0)&=&-\dfrac{a+b^2-2 b+2}{(a+1) (b-1)^2},\\
\det(H_1)&=&\dfrac{-a+4 b-1}{(a+1) (b-1)^2},\\
\det(H_2)&=&\frac{16}{(a+1) (b-1)^2}\ne 0.
\end{eqnarray*}
Therefore if $a+b^2-2 b+2\ne 0$, and $-a+4 b-1\ne 0$,  
we obtain the polynomials 
\begin{eqnarray*}
p_1(x)&=&x-\frac{a b-b^2+3 b-1}{a+b^2-2 b+2},\\
p_2(x)&=&x^2-\frac{2(a-2 b+3)}{a-4 b+1}x
+\frac{a+8 b-3}{a-4 b+1},
\end{eqnarray*}
where $p_2(x)$ is proper with $p(x)$. Moreover,  observe 
that when s$a=3$, $b=0$ we have $p_2(x)=q(x)$.
\end{example}
Theorem \ref{thm:3} has some nice consequences, 
for example, by construction, as we pointed out in 
example \ref{ex:5}, the polynomial $p_{n-1}(x)$ in such 
construction is proper with the given $p(x)$; moreover
the following result connects our problem to the 
Gaussian quadrature formulae.
\begin{remark} \label{rem:7}
Observe that, in $\mathbb C$, if $x_1=x_2=\cdots
=x_n=a$ then 
\[
\int_{\Gamma} \frac{f(z)}{p_n(z)p_{n-1}(z)}\, dz=
\frac {2\pi i}{(n-1)!} \frac{d^{n-1}}{z^{n-1}} 
\left(\frac {f(z)} {p_{n-1}(z)} \right)(a),
\]
where $\Gamma$ is a Jordan curve such that $a$ lies inside 
$\Gamma$, and the roots of $p_{n-1}$ lie outside of $\Gamma$.
\end{remark}
\begin{remark} \label{rem:9}
Taking into account theorem \ref{thm:3} and 
the Remark \ref{rem:7} if we have a field $\mathbb F$ 
in which the derivative may not make sense,
for example if $n=2$, then \eqref{3:7} becomes
\[
{\mathscr L}(f(t))=C \left(\frac {f'(a)}{p_1(a)}-\frac{f(a)}{p_1^2(a)} \right),
\]
understanding that $f(a)$ (resp. $f'(a)$) represents the 
coefficient of $(z-a)^0\equiv 1$ (resp. $(z-a)$) in the 
expansion of $f(z)$ in terms of $\{(z-a)^k\}_{k=0}^\infty$ 
in $\mathbb F$.

We can proceed In an analogous way for the the $n=3$ case. 
In such a case we have
\[
{\mathscr L}(f(t))=C \left(\frac {f''(a)}{p_2(a)}-2\frac{f'(a)p_2'(a)}
{p_2^2(a)} -2\frac{f(a)}{p_2^2(a)}+2\frac{f(a)(p_2')^2(a)}
{p_2^3(a)}\right).
\]
\end{remark}
\section{Further observations}
In this section we present some other results that are 
connected with the results presented previously.

\begin{thm1} \label{thm:10}
Let $p$ and $q$ be monic polynomials, of degree $n$ 
and $n-1$ respectively, over a field $\mathbb F$. 

If there exists an irreducible tridiagonal matrix $A$ 
such that $p_A(t)=p$ and $p_{A(1)}(t)=q$, i.e. $p$ and $q$ are 
ITrP, then
\begin{enumerate}
\item[a)]  $S_2(A)-S_2(A(1))-a_{1}S_1(A(1))\ne 0$.
\item[b)] For $k=2, 3, \dots, n-2$,
\[
S_2(A(1,\dots,k-1))-S_2(A(1,\dots,k))
-a_{k} S_1(A(1,\dots,k-1))\ne 0,
\]
\end{enumerate}
where 
\[
\det(A-\lambda I)=\sum_{k=0}^n (-1)^k S_{n-k}(A) \lambda^{k}.
\]
Conversely, if conditions a) and b) holds, then 
$p_A(t)$ and $p_{A(1)}(t)$ is a proper pair.
Note that $S_0(A)=1$, $S_1(A)=$Tr($A$), 
$S_n(A)=\det(A)$, $a_{1}=S_1(A)-S_1(A(1))$, and 
for $k=2,\dots, n-2$,
\[
a_{k}=S_1(A(1,\dots,k-1))-S_1(A(1,\dots,k)).
\]
\end{thm1}

In fact, the given conditions in theorem \ref{thm:10} 
{\it b)} can be expressed in terms of the coefficients 
of the characteristic polynomial of $A$ and $A(1)$. 
For example,  if $n=4$ such condition for $k=2$ 
can be written as follows:
\[\begin{split}
\big(S_3(A)-S_3(A(1))-  a_1 S_2(A(&1))\big)^2+b_1 
\big(S_3(A)-S_3(A(1))-a_1 S_2(A(1))\big)\\ 
&+b_1 S_2(A(1))-b_1 \big(a_1 S_3(A(1))-S_4(A)\big)\ne 0,
\end{split}\]
where $b_1=S_2(A(1))-S_2(A)+a_{1}S_1(A(1))$.
\begin{proof} 
This result follows straightforwardly by using the fact 
that $b_k$'s in the matrix $A$ can be computed as 
the coefficient of $x^{n-k-1}$ in the polynomial 
$\det(A(1, 2, \dots,k-1)-t I)-(t-a_k)\det(A(1, 2, \dots,k-1,k)-t I)$,
and that they need to be nonzero.
\end{proof}

Note that if 
\[
P_{A(1,2,\dots,k-1)}(t)=t^{n-k+1}-c_{k1}t^{n-k}+c_{k2}t^{n-k-1}+\cdots,
\]
and
\[
P_{A(1,2,\dots,k)}(t)=t^{n-k}-d_{k1}t^{n-k-1}+d_{k2}t^{n-k-2}+\cdots,
\]
then the previous result can be written as follows.
\begin{thm1} \label{thm:15}
Let $p_n$ and $p_{n-1}$ be relatively prime monic polynomials 
over a field $\mathbb F$ of degree $n$ and $n-1$. 
Then $p_n$, $p_{n-1}$ is a proper pair if and only if 
\[
(d_{k2}-c_{k2})+(c_{k1}-d_{k1})d_{k1}\ne 0,\quad k=1, 2, \dots, n-1.
\]
In this event, $p_n$, $p_{n-1}$ is an ITrP and the realizing normalized 
tridiagonal matrix is unique.
\end{thm1}
\begin{remark} Note that in the relatively prime, proper case there is 
an iterative algorithm to construct the realizing tridiagonal matrix 
computationally.
\end{remark} 

The proof of theorem \ref{thm:15} follows from the fact that 
\[
b_{k}=(d_{k2}-c_{k2})+(c_{k1}-d_{k1})d_{k1} \quad {\rm for} \ k=1, 2, \dots n-1.
\]


Another interesting fact related with our problem is the 
following. We can find the values of the Hankel determinants 
for linear functionals in terms of the roots of $p$.
Here we present the $n=4$ case:

\begin{lemma} \label{lem:13} Let $a$, $b$, $c$  and $d$ be four 
different numbers, and let $\omega_a$, $\omega_b$, $\omega_c$, 
$\omega_d$ be another four nonzero numbers. Then the 
Hankel determinants associated with the linear functional 
\[
{\mathscr L}(f)=\omega_a f(a)+\omega_b f(b)+\omega_c f(c)+
\omega_d f(d),
\]
are
\[\begin{array}{r@{\hspace{0.6mm}}c@{\hspace{0.6mm}}l}
\det(H_0)&=&\omega_a+\omega_b+\omega_c+\omega_d, \\[1mm] 
\det(H_1)&=&\omega_a \omega_b(b-a)^2+\omega_a\omega_c(c-a)^2
+\omega_a\omega_d(d-a)^2+\omega_b\omega_c(c-b)^2 \\
&& +\omega_b\omega_d(d-b)^2+\omega_c\omega_d(d-c)^2,  \\[1mm] 
\det(H_2)&=&\omega_a \omega_b\omega_c(b-a)^2(c-a)^2(c-b)^2
+\omega_a \omega_b\omega_d(b-a)^2(d-a)^2(d-b)^2 \\
&+&\omega_a\omega_c\omega_d(c-a)^2(d-a)^2(d-c)^2
+\omega_b\omega_c\omega_d(c-b)^2(d-b)^2(d-c)^2,  \\[1mm] 
\det(H_3)&=&\omega_a \omega_b\omega_c
\omega_d(b-a)^2(c-a)^2(d-a)^2(c-b)^2(d-b)^2(d-c)^2, \\[1mm] 
\det(H_k)&=& 0,\qquad k= 4, 5, \dots
\end{array}\]
\end{lemma} 

In fact, we consider the following conjecture for such values for 
the determinant of the Hankel matrices.

{\bf Conjecture:} For any different numbers 
$x_1$, $x_2, \dots, x_n$ and for any  $\omega_1, \cdots, 
\omega_n$, all different from zero, let us consider the 
linear functional 
\[
{\mathscr L}(f)=\sum_{i=1}^n \omega_i f(x_i).
\]
Then, the determinant of the Hankel matrices associated 
with this linear functional can be computed explicitly as
\[
\det (H_k)=\sum_{\substack{\Omega\subseteq \{1,2,\dots,n\}\\ 
|\Omega|=k+1}} S_{k+1}(\omega_\Omega) V^2(x_\Omega),
\]
where 
$S_k$ is the $k$-th elementary symmetric function, 
$\lambda_\Omega$ represents the set 
$\{\lambda_j: j\in \Omega\}$, and $V(\lambda_\Omega)$ 
represents the Vandermonde determinant associated with 
the numbers of the set $\lambda_\Omega$.

\begin{remark} 
Note that if $k+1>n$ then $\det (H_k)=0$, and $V(\{x\})=1$.  
\end{remark}
%
%
%

We also have considered some cases in which 
$p$ has a multiple zero. 

\begin{lemma} \label{lem:14}
For any given number $a$ and any nonzero value 
$\omega$,  let us define the linear functional (see Remark 
\ref{rem:9})
\[
{\mathscr L} (f)=\omega f'(a)-\omega^2 f(a).
\]
Then the first moments are $m_0=-\omega^2$, 
$m_1=\omega-\omega^2 a$, 
$m_2=2\omega a-\omega^2 a^2$, and
\[
\det(H_0)=\det(H_1)=-\omega^2,\quad 
\det(H_k)=0,\quad k=2, 3, \dots.
\]
\end{lemma} 

\begin{lemma} \label{lem:15}
For any given number $a$ and any two nonzero 
values $\omega_1$, $\omega_2$,  let us define the linear 
form (see Remark \ref{rem:9})
\[
{\mathscr L} (f)=\omega_1 f''(a)-2\omega_2\omega_1^2 f'(a)
-2\omega^2_1 f(a)+2\omega^2_2\omega^3_1 f(a).
\]
Then the first moments are $m_0=-2\omega_1^2+2\omega_1^3\omega_2^2$, 
$m_1=2 a \omega_2^2 \omega_1^3-2 a \omega_1^2-2 \omega_2 
\omega_1^2$, 
$m_2=2 a^2 \omega_2^2 \omega_1^3-2 a^2 \omega_1^2-4 a 
\omega_2 \omega_1^2+2 \omega_1$, and
\begin{eqnarray*}
\det(H_0)&=&2 \omega_1^2 \left(\omega_1 \omega_2^2-1\right), \\
\det(H_1)&=&-4\omega_1^2,\\ 
\det(H_2)&=&-8\omega_1^3,\\
\det(H_k)&=&0,\quad k=3, 4, \dots.
\end{eqnarray*}
\end{lemma} 

\section*{Acknowledgements}
The first author wants to thank Prf. Charles Johnson for his hospitality 
and for the opportunity to collaborate with him in the R.E.U. 
program at the College of William and Mary. He also acknowledges 
financial support by Direcci\'on General de Investigaci\'on, Ministerio 
de Econom\'ia y Competitividad of Spain, grant MTM2015-65888-C4-2-P.
Both authors wish to thank students Colin Walker and Owen Hill for past 
work on this problem.

\end{document}